\documentclass{article}
\usepackage{amssymb}
\usepackage{amsfonts}
\usepackage{mathrsfs}
\usepackage{amsmath}
\usepackage{latexsym}
\usepackage{float}
\usepackage{color}
\usepackage[utf8]{inputenc}
\usepackage[pdftex]{graphicx}
\newtheorem{theorem}{Theorem}
\newtheorem{lemma}{Lemma}

\newcommand{\R}{{\mathbb{R}}}
\newcommand{\C}{{\mathbb{C}}}

\newcommand{\cL}{{\cal L}}

\newcommand{\cE}{{\cal E}}

\newcommand{\po}{\partial\Omega}
\newcommand{\re}{{\rm Re}\;}

\newcommand{\tarr}[6]{{\left\{\begin{array}{lll}
{#1}&{#2}\\[0.2cm]
{#3}&{#4}\\[0.2cm]
{#5}&{#6}
\end{array}\right.}}

\newcommand{\ia}{({\rm i})}
\newcommand{\ib}{({\rm ii})}
\newcommand{\ic}{({\rm iii})}
\newcommand{\id}{({\rm iv})}

\newcommand{\lf}{\lambda\varphi}

\newcommand{\be}{\begin{equation}}
\newcommand{\ee}{\end{equation}}
\newcommand{\bea}{\begin{eqnarray}}
\newcommand{\eea}{\end{eqnarray}}
\newcommand{\bean}{\begin{eqnarray*}}
\newcommand{\eean}{\end{eqnarray*}}
\newcommand{\la}{\label}

\title{On the heat kernel of a class of fourth order operators in two dimensions:
sharp Gaussian estimates and short time asymptotics}
\date{}
\author{
G. Barbatis\thanks{Department of Mathematics,
 National and Kapodistrian University of Athens,  15784 Athens, Greece}
 \and P. Branikas\footnotemark[1]
}
\begin{document}
\date{}
\maketitle
\vspace{11pt}
\parindent=0pt
\parskip=8pt

\begin{abstract}
We consider a class of fourth order uniformly elliptic operators in planar Euclidean domains
and study the associated heat kernel. For operators with $L^{\infty}$ coefficients we obtain Gaussian estimates with best constants, while for operators with constant coefficients we obtain short time asymptotic estimates. The novelty of this work is that we do not assume that the associated symbol is strongly convex. The short time asymptotics reveal a behavior which is qualitatively different from that of the strongly convex case. 
\end{abstract}

\vspace{11pt}


%\vspace{6pt}

\section{Introduction}

Let $\Omega$ be a planar domain and let
\[
Hu = \partial_{x_1}^2 \big( \alpha(x) \partial_{x_1}^2 u  \big) +2\partial_{x_1}\partial_{x_2} \big( \beta(x) \partial_{x_1}\partial_{x_2} u  \big)
+ \partial_{x_2}^2 \big( \gamma(x) \partial_{x_2}^2 u  \big)
\]
be a self-adjoint, fourth-order uniformly elliptic operator in divergence form on $\Omega$ with $L^{\infty}$ coefficients
satisfying Dirichlet boundary conditions on $\partial\Omega$. It has been proved by Davies \cite{davies1} that $H$ has a heat kernel
$G(x,x',t)$ which satisfies the Gaussian-type estimate,
\be
|G(x,x',t)|\leq c_1 t^{-\frac{1}{2}}\exp\Big(-c_2\frac{|x-x'|^{4/3}}{t^{1/3}}+c_3t\Big),
\la{in0}
\ee
for some positive constants $c_1$, $c_2$, $c_3$ and all $t>0$ and $x,y\in\Omega$.

The problem of finding the sharp value of the exponential constant $c_2$ is related to replacing the Euclidean distance $|x-y|$ by an appropriate distance $d(x,y)$ that is suitably adapted to the operator $H$ and, more preciesly, to its symbol
\[
A(x,\xi) =\alpha(x) \xi_1^4 +2\beta(x) \xi_1^2 \xi_2^2 +\gamma(x)\xi_2^4 \; , \qquad x\in\Omega \,  ,  \;\; \xi\in\R^2 \, .
\]
In the article \cite{ep}, and for constant coefficient operators in $\R^n$ which satisfy suitable assumptions, the asymptotic formula
\be
G(x,x',t)\sim h(x-x')^{-2/3} t^{-1/3}\exp\Big(- \frac{3\sqrt[3]{2}}{16} \frac{p_*(x-x')^{4/3}}{t^{1/3}}\Big)
\cos\Big(-\frac{3\sqrt{3}\sqrt[3]{2}}{16} \frac{p_*(x-x')^{4/3}}{t^{1/3}} -\frac{\pi}{3}\Big), 
\la{in3}
\ee
was established as $t\to 0+$; here $h$ is a positively homogeneous function of degree one and $p_*$ is the Finsler metric defined by
\be
p_*(\xi) = \max_{\eta \in\R^2\setminus\{0\}} \frac{ \eta\cdot \xi}{ A(\eta)^{1/4}}.
\la{plo}
\ee
An analogous asymptotic formula has been obtained in \cite{tintarev} in the more general case of operators with variable smooth coefficients;
in this case the relevant distance is the (geodesic) Finsler distance $d_{p_*}(x,x')$ induced by the Finsler metric with length element $p_*(x,\xi)$, the latter being defined similarly to (\ref{plo}), with the additional dependence on $x$.

A sharp version of the Gaussian estimate (\ref{in0}) was established in \cite{b2001} where it was proved that
\be
|G(x,x',t)|\leq c_{\epsilon} t^{-\frac{1}{2}}\exp\Big\{-\big( \frac{3\sqrt[3]{2}}{16} -D -\epsilon\big) \frac{d_M(x,x')^{4/3}}{t^{1/3}}+
c_{\epsilon,M}t\Big\},
\la{in2}
\ee
for arbitrary $\epsilon$ and $M$ positive. Here $D\geq 0$ is a constant that is related to the regularity of the coefficients and $d_M(x,x')$, $M>0$, is a family of Finsler-type distances on $\Omega$ which
is monotone increasing and converges as $M\to +\infty$ to a limit Finsler-type distance $d(x,x')$ closely related to $d_{p_*}(x,x')$ but not equal to it; see also Subsection \ref{nss}.

A fundamental assumption for both (\ref{in3})  and (\ref{in2})  is the {\em strong convexity} of the symbol $A(x,\xi)$ of the operator $H$.
The notion of strong convexity was introduced in \cite{ep} where short time asymptotics were obtained not only for the operator described above but more generally for a constant coefficient operator of order $2m$ acting on functions on $\R^n$.
In the context of the present article and for an operator with constant coefficients, strong convexity of the symbol $A(\xi)$ amounts to
\be
\la{sc1}
0< \beta  < 3 \sqrt{\alpha\gamma} \, .
\ee
We note however that in \cite{b2001} (where the coefficients $\alpha,\beta,\gamma$ are functions) the term strong convexity was also used for the slightly more general situation
where
\be
\la{sc2}
0\leq  \beta(x) \leq 3 \sqrt{\alpha(x)\gamma(x)} \, , \qquad x\in \Omega. 
\ee
In other words, while for short time asymptotics the strict inequality was necessary, for Gaussian estimates equality is allowed.

Our aim in this work is to extend both (\ref{in3}) and (\ref{in2}) to the case of non-strongly convex symbols.
Hence in Theorem 1, which extends the Gaussian estimate (\ref{in2}),
condition (\ref{sc2}) is not valid, while in Theorem 2, which extends the short time asymptotics (\ref{in3}), condition (\ref{sc1}) is not valid.

In Theorem \ref{thm1} we obtain a Gaussian estimate involving the Finsler-type distances $d_M(\cdot , \cdot)$ and an  $\sigma_*$ that depends on the range of the function
\[
Q(x)=\frac{\beta(x)}{\sqrt{\alpha(x)\gamma(x)}} \, .
\]
So the strongly convex case corresponds to $Q$ taking values in $[0,3]$ but here we allow $Q$ to take any value in $(-1, +\infty)$. It is worth noting that while in the strongly convex case we have $\sigma_*  =3\sqrt[3]{2} /16$, in the general case
$\sigma_*$ can take a continuous range of values.
Our approach follows the main ideas of \cite{b2001} and in particular makes use of Davies' exponential perturbation method.
However technical difficulties arise due to the existence of three different regimes for the function $Q(x)$,
namely $(-1,0]$, $[0,3]$ and $[3,+\infty)$. Each regime must be handled differently,
and it must be shown that the matching at $Q=0$ and $Q=3$ does not cause any problems.

In the second part of the paper we extend the short time asymptotic estimates of \cite{ep} to operators with non-strongly convex symbol.
As in \cite{ep}, we consider constant coefficient operators
acting on $\R^2$, so the heat kernel (also referred to as the Green's function) is given by
\[
G(x,t)=\frac{1}{(2\pi)^2}\int_{\R^2}e^{i\,\xi\cdot\,x -A(\xi) t}\,d\xi,\qquad x\in\R^2\, ,\; t>0 .
\]
The asymptotic estimates are contained in Theorem \ref{thm2} and the proof uses the steepest descent method.
For technical reasons we only consider specific choices of the point $x\in\R^2$; we comment further upon this before the statement of Theorem \ref{thm2},
but note that the asymptotic formulae obtained are enough to demonstrate the sharpness of the exponential constant of Theorem \ref{thm1}.
Anyway, these asymptotic estimates are of independent interest, in particular because they reveal a behavior that is qualitatively different from that of the strongly convex case.
In the case $0< Q<3$ studied in \cite{ep} the Green function oscillates around the horizontal axis. As it turns out, when
$ Q<0$ or $ Q>3$ the Green function remains positive for small times. The borderline
cases $ Q=0$ and $ Q=3$ are particularly interesting. In these two cases the asymptotic expression
involves oscillations that touch the horizontal axis at their lowest points (see the diagrams at the end of the article). This is due to a bifurcation phenomenon that takes place at $Q=0$ and $Q=3$.
At these values of $ Q$ there is a change in the branches of saddle points that contribute to the asymptotic behavior of the integral. While for each $ Q\neq 0,3$ there are two contributing points,
for each of the values $ Q=0$ and $ Q= 3$ there are four such points.

At the end of the article we present numerical computations that illustrate the asymptotic estimates. For the sake of completeness we have also included an appendix with the proof of Evgrafov and Postnikov in the strongly convex case $0<Q<3$.

We close this introduction with one remark. As pointed out in \cite{ep}, for a fourth order operator in two dimensions
strong convexity is equivalent to convexity. Nevertheless we have chosen not to replace the term `strongly convex' by `convex' in order to
emphasize the importance of strong convexity in the general case of an operator of order $2m$ acting in $\R^n$ (considered in both \cite{ep} and \cite{b2001}).

\section{Heat kernel estimates}

\subsection{Setting and statement of Theorem \ref{thm1}}

Let $\Omega\subset\R^2$ be open and connected. We consider a differential operator $H$ on $L^2(\Omega)$ (complex-valued functions) given formally by
\[
Hu(x)=\partial_{x_1}^2\big(\alpha(x)\partial_{x_1}^2u\big)+2\partial^2_{x_1x_2}(\beta(x)\partial^2_{x_1x_2}u)+\partial_{x_2}^2\big(\gamma(x)\partial_{x_2}^2u\big),
\]
where $\alpha$, $\beta$ and $\gamma$ are functions in $L^{\infty}(\Omega)$. In case $\Omega\neq \R^2$ we impose
Dirichlet boundary conditions on $\po$. 
The operator $H$ is defined by means of the quadratic form
\[
Q(u)=\int_{\Omega}\big\{\alpha(x)|u_{x_1x_1}|^2+2\beta(x)|u_{x_1x_2}|^2+\gamma(x)|u_{x_2x_2}|^2\big\}\,dx,
\]
defined on $H^2_0(\Omega)$. We assume that $H$ is uniformly elliptic, that is the functions $\alpha$ and $\gamma$ are positive and bounded away from zero and also
\[
\inf_{x\in\Omega} \frac{\beta(x)}{\sqrt{\alpha(x)\gamma(x)}} > -1 \, .
\]
The form $Q$ is then closed and we define the operator $H$ to be the self-adjoint operator associated to it. As mentioned in the introduction, the operator $H$ has a heat kernel $G(x,x',t)$ which satisfies the Gaussian estimate (\ref{in0}).

To state the main result of this section we need to introduce some more definitions.
We define the class of real-valued functions
\[
\cE =\{\varphi\in C^2(\Omega):\;\|D^\alpha\varphi\|_\infty<+\infty\,,\;0\leq|\alpha|\leq2\},
\]
and the subclass
\[
\cE_{A,M} =\{\varphi\in\cE \, :\;A(y,\nabla\varphi(y))\leq1\, , \;
y\in\Omega \,  , \; \mbox{ and }\;\|D^\alpha\varphi\|_{\infty}\leq M\,,\;|\alpha|=2\}.
\]
We then define a distance $d_M(\cdot,\cdot)$ on $\Omega$ by
\[
d_M(x,x')=\sup\big\{\varphi(x')-\varphi(x)\,:\;\;\varphi\in\cE_{A,M}\big\}.
\]
It is not difficult to see that as $M\to +\infty$ this converges to the Finsler-type distance
\be
d(x,x') = \sup \{ \varphi(x') -\varphi(x) \; : \;  \varphi\in{\rm Lip}(\Omega) , \;\; A( y ,\nabla\varphi(y))\leq 1 \, , \;\; y\in\Omega\}.
\la{fd}
\ee
As it turns out, there holds $d(x,x')\leq d_{p_*}(x,x')$ and the two distances in general are not equal. Still, there are cases where equality is valid and this shows in particular that the best constant for Gaussian estimates is the same for both distances.
This is further discussed in Subsection \ref{nss}.

We next define the functions
\[
Q(x)=\frac{\beta(x)}{\sqrt{\alpha(x)\gamma(x)}},\qquad\quad x\in\Omega,
\]
\[
k(x)=\tarr{8\frac{1-Q(x)}{(1+Q(x))^2},}{\mbox{ if }-1<Q(x)<0,}{8,}{\mbox{ if }0\leq Q(x)\leq3,}{Q(x)^2-1,}{\mbox{ if }Q(x)>3,}
\]
and
\be
\sigma(x)=\frac{3}{4}\cdot\Big(\frac{1}{4k(x)}\Big)^{1/3}= \tarr{ \frac{3}{2\cdot 4^{4/3}} \frac{ (1+Q(x))^{2/3}}{ (1-Q(x))^{1/3}},} 
{\mbox{ if }-1<Q(x)<0,}{\frac{3}{8\cdot 4^{1/3}},}{\mbox{ if }0\leq Q(x)\leq3,}{ \frac{3}{4^{4/3}} (Q(x)^2-1)^{-1/3} ,}{\mbox{ if }Q(x)>3,}
\la{sx}
\ee
We set 
\[
k^*=\sup_{x\in\Omega}k(x)  \, , \quad \mbox{ and }
\sigma_*= \inf_{x\in\Omega}\sigma(x) = \frac{3}{4}\cdot\Big(\frac{1}{4k^*}\Big)^{1/3}.
\]
Finally, we denote by $D$ the distance in $L^{\infty}(\Omega)$ of the functions $\alpha(x),\beta(x),\gamma(x)$ from the space of all Lipschitz functions,
\be
D=\max\big\{d_{L^{\infty}}(\alpha,{\rm Lip}(\Omega))\,,\;d_{L^{\infty}}(\beta,{\rm Lip}(\Omega))\,,\;d_{L^{\infty}}(\gamma,{\rm Lip}(\Omega))\big\}.
\la{D}
\ee
The main result of this section is the following
\begin{theorem}
For all $\epsilon\in(0,1)$ and all $M$ large there exists $c_\epsilon,c_{\epsilon,M}<\infty$ such that
\be
|G(x,x',t)|\leq c_\epsilon t^{-1/2}\exp\Big\{-(\sigma_*-cD-\epsilon)\frac{d_M(x,x')^{4/3}}{t^{1/3}}+c_{\epsilon,M}t\Big\},
\la{eq:mainthm}
\ee
for all $x,x'\in\Omega$ and $t>0$.
\la{thm1}
\end{theorem}
It will follow from the results of Section \ref{sec:asympt} that the constant $\sigma_*$ is sharp.

%%%%%%%%%%%%%%%%%%%%%%%%%%%%%
%%%%%%%%%%%%%%%%%%%%%%%%%%%%%

\subsection{Proof of Theorem \ref{thm1}}

We first establish some auxiliary inequalities related to the symbol $A(x,\xi)$ of the operator $H$. Since these are pointwise inequalities with respect to $x\in\Omega$, we assume for simplicity that the coefficients are constant and therefore the symbol is
\[
A(\xi)=\alpha\xi_1^4+2\beta\xi_1^2\xi_2^2+\gamma\xi_2^4,\qquad\xi\in\R^2,
\]
where $\alpha,\beta,\gamma\in\R$. By ellipticity we have $\alpha,\gamma>0$ and $Q:=\beta/\sqrt{\alpha\gamma}>-1$.
We shall need to consider the symbol also as a function of two complex variables, that is we set
\[
A(z)=\alpha z_1^4+2\beta z_1^2 z_2^2+\gamma z_2^4,\qquad z=(z_1,z_2)\in\C^2 \, .
\]
\begin{lemma}
There holds
\be
\la{eq1}
\re\,A(\xi+i\eta)\geq-kA(\eta),\quad\quad \xi,\,\eta\in\R^2,
\ee
where the constant $k$ is given by
\[
k=\tarr{8\frac{1-Q}{(1+Q)^2},}{\mbox{ if }\;-1<Q<0,}{8,}{\mbox{ if }\;0\leq Q\leq 3,}{Q^2-1,}{\mbox{ if }\;Q>3.}
\]
\label{lem1}
\end{lemma}
{\em Proof.} We first compute
\bea
\re\,A(\xi+i\eta)&=& \alpha (\xi_1^4 -6\xi_1^2\eta_1^2 +\eta_1^4 ) +2\beta\Big(\xi_1^2\xi_2^2
-\xi_1^2\eta_2^2-\xi_2^2\eta_1^2-4\xi_1\xi_2\eta_1\eta_2+\eta_1^2\eta_2^2\Big) \nonumber\\
&& + \gamma (\xi_2^4   -6\xi_2^2\eta_2^2+\eta_2^4 ).
\la{expansion}
\eea

We now distinguish the three cases.

(i) $-1<Q<0$. Using (\ref{expansion}) we see by a direct computation that
\bea
&&\re\,A(\xi+i\eta)+8\frac{1-Q}{(1+Q)^2}A(\eta)\nonumber\\
&=&(Q+1)\Big\{\alpha\Big(\xi_1^2-\frac{3-Q}{1+Q}\eta_1^2 \Big)^2+\gamma\Big(\xi_2^2-\frac{3-Q}{1+Q}\eta_2^2\Big)^2\Big\}
-Q\big(\alpha^{1/2}\xi_1^2-\gamma^{1/2}\xi_2^2\big)^2 \la{relation1}\\
&&-2Q\big(\alpha^{1/2}\xi_1\eta_1+\gamma^{1/2}\xi_2\eta_2\big)^2-2Q\alpha^{1/2}\gamma^{1/2}(\xi_1\eta_2+\xi_2\eta_1)^2-Q\Big(\frac{3-Q}{1+Q}\Big)^2 \big(\alpha^{1/2}\eta_1^2-\gamma^{1/2}\eta_2^2\big)^2 \! ,\nonumber
\eea
and (\ref{eq1}) follows.

(ii) $0\leq Q\leq3$. Similarly it may be verified that
\bea
\re\,A(\xi+i\eta)+8A(\eta)&=&\frac{Q}{3}\Big\{\alpha^{1/2}(\xi_1^2-3\eta_1^2)+\gamma^{1/2}(\xi_2^2-3\eta_2^2)\Big\}^2+\frac{4Q}{3}\alpha^{1/2}\gamma^{1/2}(\xi_1\xi_2-3\eta_1\eta_2)^2\nonumber\\
&&+\frac{3-Q}{3}\Big\{\alpha(\xi_1^2-3\eta_1^2)^2+\gamma(\xi_2^2-3\eta_2^2)^2\Big\},\la{relation2}
\eea
and (\ref{eq1}) again follows.

(iii) $Q>3$. In this case we have
\bea
&&\hspace{-1.5cm}\re\,A(\xi+i\eta)+(Q^2-1)A(\eta)\nonumber\\
&=&2(Q-3)\big(\alpha^{1/2}\xi_1\eta_1-\gamma^{1/2}\xi_2\eta_2\big)^2  +\Big\{\alpha^{1/2}(\xi_1^2-Q\eta_1^2)+\gamma^{1/2}(\xi_2^2-Q\eta_2^2)\Big\}^2 \la{relation3}\\
&&+2(Q-1)\alpha^{1/2}\gamma^{1/2}\Big(\xi_1\xi_2-\frac{Q+3}{Q-1}\eta_1\eta_2\Big)^2+2\frac{(Q-3)(Q+1)(Q^2+3)}{Q -1}\alpha^{1/2}\gamma^{1/2}\eta_1^2\eta_2^2, \nonumber
\eea
and (\ref{eq1}) follows once again.  $\hfill\Box$

%%%%%%%%%%%%%%%%%%%%%%

Given $\psi\in\cE$ the (multiplication) operator $e^{\psi}$ leaves the Sobolev space $H_0^2(\Omega)$ invariant so one can define a quadratic form $Q_{\psi}$ on $H^2_0(\Omega)$ by
\[
Q_{\psi}(u)=Q(e^{\psi}u,e^{-\psi}u),
\]
where
\[
Q(u,v)= \int_{\Omega}\big\{\alpha(x)u_{x_1x_1}\overline{v}_{x_1x_1}+2\beta(x)u_{x_1x_2}\overline{v}_{x_1x_2}
+\gamma(x)u_{x_2x_2}\overline{v}_{x_2x_2}\big\}\,dx
\]
is the sesquilinear form associated to $Q(\cdot)$.
Expanding the various terms of $Q_{\psi}(u)$ (cf. (\ref{pat}) below) we find that the highest order terms
coincide with those of $Q(u)$ and standard interpolation inequalities (cf. \cite[Lemma 2]{davies1}) then give
\be
|Q(u)-Q_{\psi}(u)|\leq\epsilon Q(u)+c_{\epsilon}\{\|\psi\|_{W^{2,\infty}}+\|\psi\|_{W^{2,\infty}}^{4}\}\|u\|_2^2, 
\la{fil}
\ee
for all $\epsilon\in(0,1)$ and $u\in H^2_0(\Omega)$.

The proof of Theorem \ref{thm1} makes essential use of the following result which is Proposition 2 of \cite{bd}:
\begin{lemma}
Let $\psi\in\cE $ be given and let $\tilde k >0$ be such that
\[
{\rm Re}\,Q_{\psi}(u)\geq - \tilde{k}  \,\|u\|_2^2,
\]
for $u\in C_c^{\infty}(\Omega)$. Then for any $\delta\in(0,1)$ there exists a constant $c_\delta$ such that
\[
|G(x,y,t)|\leq c_{\delta} t^{-1/2}\exp\{\psi(x)-\psi(y)+(1+\delta) \tilde{k} t\},
\]
for all $x,y\in\Omega$ and all $t>0$.
\label{lem:jde}
\end{lemma}
Given $\varphi\in\cE$ and $\lambda>0$ we have
\bea
Q_{\lf}(u)&=&\int_{\Omega}\Big[\alpha(x)(e^{\lf}u)_{x_1x_1}
(e^{-\lf}\overline{u})_{x_1x_1}+2\beta(x)(e^{\lf}u)_{x_1x_2}(e^{-\lf}\overline{u})_{x_1x_2}\nonumber\\
&&\quad\quad+\gamma(x)(e^{\lf}u)_{x_2x_2}(e^{-\lf}\overline{u})_{x_2x_2} \Big] \,dx.
\la{pat}
\eea

Using Leibniz's rule to expand the second partial derivatives the exponentials $e^{\lf}$ and $e^{-\lf}$ cancel and we conclude that $Q_{\lf}(u)$ is a linear combination of terms of the form
\be
\la{terms}
\lambda^s\int_{\Omega}b_{s\gamma\delta}(x)D^{\gamma}u\,D^{\delta}\overline{u}\,dx,
\ee
(multi-index notation) where each function $b_{s\gamma\delta}(x)$ is a product of one of the functions $\alpha(x)$, $\beta(x)$, $\gamma(x)$ and first or second order derivatives of $\varphi$. For any such term we have $s+|\gamma+\delta|\leq4$.

{\bf Definition.} We denote by $\cL$ the space of (finite) linear combinations of terms of the form (\ref{terms})
with $s+|\gamma+\delta|  <4$.

We shall see later the terms in $\cL$ are in a certain sense negligible.
We next define the quadratic form
\bean
Q_{1,\lf}(u)\!&=&\int_{\Omega}\Big\{\lambda^4\big[\alpha(x)\varphi_{x_1}^4+2\beta(x)\varphi_{x_1}^2\varphi_{x_2}^2+\gamma(x)\varphi_{x_2}^4\big]|u|^2\\
&&\quad+\lambda^2\Big\{\alpha(x)\varphi_{x_1}^2(u\overline{u}_{x_1x_1}+u_{x_1x_1}\overline{u} -4|u_{x_1}|^2)\\
&&\quad+2\beta(x)\big[\varphi_{x_1}\varphi_{x_2}(u\overline{u}_{x_1x_2}+u_{x_1x_2}\overline{u}-u_{x_1}\overline{u}_{x_2}-u_{x_2}\overline{u}_{x_1})-(\varphi_{x_2}^2|u_{x_1}|^2+\varphi_{x_1}^2|u_{x_2}|^2)\big]\\
&&\quad+\gamma(x)\varphi_{x_2}^2(u\overline{u}_{x_2x_2}+u_{x_2x_2}\overline{u}-4|u_{x_2}|^2)\Big\}\\
&&\quad +\alpha(x)|u_{x_1x_1}|^2+2\beta(x)|u_{x_1x_2}|^2+\gamma(x)|u_{x_2x_2}|^2\Big\}\,dx.
\eean
It can be easily seen that $Q_{1,\lf}(\cdot)$ contains precisely those terms of the form (\ref{terms}) from the expansion of $Q_{\lf}(\cdot)$ for which we have $s+|\gamma+\delta|=4$. Hence we have
\begin{lemma}
The difference $Q_{\lf}(\cdot)-Q_{1,\lf}(\cdot)$ belongs  to $\cL$.
\label{lem2}
\end{lemma}

The symbol of the operator $H$ is
\[
A(x,z)=\alpha(x)z_1^4+2\beta(x)z_1^2z_2^2+\gamma(x)z_2^4,\qquad x\in\Omega,\quad z\in\C^2\,,
\]
and the polar symbol is defined as
\[
A(x,z,z')=\alpha(x)z_1^2z_1'^2+2\beta(x)z_1z_2z_1'z_2'+\gamma(x)z_2^2z_2'^2,\qquad x\in\Omega,\quad z,\,z'\in\C^2.
\]
For $x\in\Omega$ and $\xi,\xi',\eta\in\R^2$ we set
\[
S(x,\xi,\xi',\eta)={\rm Re}\,A(x,\xi+i\eta,\xi'+i\eta)+k(x)A(x,\eta).
\]
Given $\varphi\in\cE$ and $\lambda\in\R$ we define the quadratic form $S_{\lf}$ on $H^2_0(\Omega)$ by
\[
S_{\lf}(u)=\frac{1}{(2\pi)^{2}}\iiint_{\Omega\times\R^2\times\R^2}S(x,\xi,\xi',\lambda\nabla\varphi)e^{i(\xi-\xi')\cdot x}\hat{u}(\xi)\overline{\hat{u}(\xi')}\,d\xi\,d\xi'\,dx.
\]
\begin{lemma}
There holds
\[
Q_{1,\lf}(u)+\int_{\Omega} k(x)A(x,\lambda\nabla\varphi)|u|^2dx=S_{\lf}(u),
\]
for all $\varphi\in\cE$, $\lambda>0$ and $u\in H^2_0(\Omega)$.
\la{lem:4}
\end{lemma}

{\em Proof.} For the proof one simply uses the relation $D^{\alpha}u(x) =(2\pi)^{-1}\int_{\R^2}(i\xi)^{\alpha}e^{ix\cdot\xi}\hat{u}(\xi)d\xi$
for the various terms that appear in $Q_{1,\lf}$. Since a very similar proof has been given in \cite{b2001} we omit further details (the fact that $k(x)$ is not constant in our case is not a problem and strong convexity is not relevant here). $\hfill\Box$

We now define for each $x\in\Omega$ a quadratic form $\Gamma(x , \cdot)$ in $\C^6$ by
\bean
&&\Gamma(x,p)=\\
&&\hspace{-0.5cm}=\left\{\begin{array}{lll}{(Q+1)|p_1|^2+(Q+1)|p_2|^2-Q|p_3|^2-2Q|p_4|^2-}&{}&{}\\[0.1cm]
{\hspace{3cm}-2Q|p_5|^2-\frac{Q(3-Q)^2}{(1+Q)^2}|p_6|^2,}&&{\hspace{-0.3cm}\mbox{if }-1<Q(x)<0,}\\[0.3cm]
{\frac{3-Q}{3}|p_1|^2+\frac{3-Q}{3}|p_2|^2+\frac{Q}{3}|p_1+p_2|^2+\frac{4Q}{3}|p_3|^2,}&{}&
{\hspace{-0.3cm}\mbox{if }0\leq Q(x)\leq3,}\\[0.3cm]
{2(Q-3)|p_1|^2+|p_2|^2 +2(Q-1)|p_3|^2+2\frac{Q-3}{Q-1}(Q+1)(Q^2+3)|p_4|^2,}&&{\hspace{-0.3cm}\mbox{if }Q(x)>3.}
\end{array}
\right.
\eean
for any $p=(p_1,\ldots,p_6)\in\C^6$.
Clearly $\Gamma(x , \cdot)$ is positive semidefinite for each $x\in\Omega$. We denote by $\Gamma(x,\cdot,\cdot)$ the corresponding sesquilinear form in $\C^6$, that is $\Gamma(x, p ,q)$ is given by a formula similar to the one above with each $|p_k|^2$ being replaced by $p_k\overline{q_k}$.

Next, for any $x\in\Omega$ and $\xi,\eta\in\R^2$ we define a vector $p_{x,\xi,\eta}\in\R^6$ by
\bean
&&p_{x,\xi,\eta}=\\
&& \left\{
\begin{array}{l}
{\Big(\alpha^{1/2}[\xi_1^2-\frac{3-Q}{1+Q}\eta_1^2],\,\gamma^{1/2}[\xi_2^2-\frac{3-Q}{1+Q}\eta_2^2],\,\alpha^{1/2}\xi_1^2-\gamma^{1/2}\xi_2^2,\,\alpha^{1/2}\xi_1\eta_1+\gamma^{1/2}\xi_2\eta_2,}\\
{\qquad \alpha^{1/4}\gamma^{1/4}(\xi_1\eta_2+\xi_2\eta_1),\,\alpha^{1/2}\eta_1^2-\gamma^{1/2}\eta_2^2\Big),}
\hspace{3.3cm}{\mbox{ if } -1<Q(x)<0,} \\[0.2cm]
{\Big(\alpha^{1/2}[\xi_1^2-3\eta_1^2],\,\gamma^{1/2}[\xi_2^2-3\eta_2^2],\,\alpha^{1/4}\gamma^{1/4}[\xi_1\xi_2-3\eta_1\eta_2],\,0,\,0,\,0\Big),}\hspace{0.5cm}{\mbox{ if } 0\leq Q(x)\leq3,} \\[0.2cm]
{\Big(\alpha^{1/2}\xi_1\eta_1-\gamma^{1/2}\xi_2\eta_2,\,\alpha^{1/2}(\xi_1^2-Q \eta_1^2)+\gamma^{1/2}(\xi_2^2-Q\eta_2^2),}\\
{\qquad\alpha^{1/4}\gamma^{1/4}[\xi_1\xi_2-\frac{Q+3}{Q-1}\eta_1\eta_2],\,\alpha^{1/4}\gamma^{1/4}\eta_1\eta_2,\,0,\,0\Big),}
\hspace{2.45cm} {\mbox{ if }Q(x)>3.}
\end{array}
\right.
\eean
A crucial property of the form $\Gamma(x,\cdot)$ and the vectors $p_{x,\xi,\eta}$ is that
\be
S(x;\xi,\xi,\eta)=\Gamma(x,p_{x,\xi,\eta},p_{x,\xi,\eta}),
\label{s:g}
\ee
for all $x\in\Omega$ and $\xi,\eta\in\R^2$; this is an immediate consequence of relations (\ref{relation1}), (\ref{relation2})
and (\ref{relation3}), for each of the three cases respectively.

We next define a quadratic form $\Gamma_{\lf}(\cdot)$ on $H^2_0(\Omega)$ by
\[
\Gamma_{\lf}(u)=\frac{1}{(2\pi)^{2}}\iiint_{\Omega\times\R^2\times\R^2}\Gamma(x, \, p_{x,\xi,\lambda\nabla\varphi},p_{x,\xi',\lambda\nabla\varphi})e^{i(\xi-\xi')\cdot x}\hat{u}(\xi)\overline{\hat{u}(\xi')}\,d\xi\,d\xi'\,dx.
\]
We then have
\begin{lemma}
Assume that the functions $\alpha(x),\beta(x),\gamma(x)$ are Lipschitz continuous. Then the difference $S_{\lf}(\cdot)-\Gamma_{\lf}(\cdot)$ belongs to $\cL$.
\label{lem3}
\end{lemma}
{\em Proof.} We consider the difference
\[
S(x,\xi,\xi',\eta)-\Gamma(x,p_{x,\xi,\eta},p_{x,\xi',\eta}),
\]
and we group together terms that have the property that if we set $\xi'=\xi$ then they are similar as monomials of the variables
$\xi$ and $\eta$. Due to (\ref{s:g}) one can use integration by parts to conclude that the total contribution of each such
group belongs to $\cL$. We shall illustrate this for one particular group, the one consisting of terms
which for $\xi=\xi'$ involve the term  $\xi_1^2\eta_1^2$.

The terms of this group from $S(x,\xi,\xi',\eta)$ add up to
\[
-\alpha(x)\eta_1^2(\xi_1^2+\xi_1'^2+4\xi_1\xi_1').
\] 
The corresponding terms of $\Gamma(x,p_{x,\xi,\eta},p_{x,\xi',\eta})$ are
\[
\tarr{\alpha(x)\eta_1^2\big[(Q(x)-3)(\xi_1^2+\xi_1'^2)-2Q(x)\xi_1\xi_1'\big],}{\mbox{ if }Q(x)<0,}
{-3\alpha(x)\eta_1^2(\xi_1^2+\xi_1'^2),}{\mbox{ if }0\leq Q(x)\leq3,}
{\alpha(x)\eta_1^2\big[-Q(x)(\xi_1^2+\xi_1'^2)+2(Q(x)-3)\xi_1\xi_1'\big],}{\mbox{ if}Q(x)>3.}
\]
Hence the difference of these terms in $S(x,\xi,\xi',\eta)-\Gamma(x,p_{x,\xi,\eta},p_{x,\xi',\eta})$ is
\[
\tarr{\alpha(x)\eta_1^2\big[(2-Q(x))(\xi_1^2+\xi_1'^2)+(2Q(x)-4)\xi_1\xi_1'\big],}{\mbox{ if }Q(x)<0,}
{\alpha(x)\eta_1^2\big[2(\xi_1^2+\xi_1'^2)-4\xi_1\xi_1'\big],}{\mbox{ if }0\leq Q(x)\leq3,}
{\alpha(x)\eta_1^2\big[(Q(x)-1)(\xi_1^2+\xi_1'^2)+(2-2Q(x))\xi_1\xi_1'\big],}{\mbox{ if } Q(x)>3.}
\]

This can also be written as $\alpha(x)\eta_1^2 R(x)(\xi_1^2+\xi_1'^2-2\xi_1\xi_1')$ where
\[
R(x)=\tarr{2-Q(x),}{\mbox{ if }Q(x)<0,}{2,}{\mbox{ if }0\leq Q(x)\leq3,}{Q(x)-1,}{\mbox{ if }Q(x)>3.}
\]

Inserting this in the triple integral and recalling that $\eta=\lambda\nabla\varphi$ we obtain that the contribution of the above terms in the difference $S_{\lf}(u)-\Gamma_{\lf}(u)$ is
\bean
&&(2\pi)^{-2}\iiint_{\Omega\times\R^2\times\R^2}\alpha(x)R(x)(\xi_1^2+\xi_1'^2-2\xi_1\xi_1')\lambda^2\varphi_{x_1}^2e^{i(\xi-\xi')\cdot x}\hat{u}(\xi)\overline{\hat{u}(\xi')}\,d\xi\,d\xi'\,dx\\
&=&\lambda^2\int_{\Omega}\alpha(x)R(x)\varphi_{x_1}^2(-u_{x_1x_1}\overline{u}-u\overline{u}_{x_1x_1}-2|u_{x_1}|^2)dx\\
&=&-\lambda^2\int_{\Omega}\alpha(x)R(x)\varphi_{x_1}^2( u_{x_1}\overline{u}+u\overline{u}_{x_1})_{x_1}dx.
\eean
Since the function $\alpha(x)R(x)\varphi_{x_1}^2$ is Lipschitz continuous we can integrate by parts and conclude that the last integral belongs to $\cL$. Similar considerations are valid for the other groups; we omit further details. $\hfill\Box$
\begin{lemma}
Assume that the functions $\alpha(x),\beta(x),\gamma(x)$ are Lipschitz continuous. Let $M>0$ be given. Then for any $\varphi\in\cE_{A,M}$ and $\lambda>0$ we have
\[
\re\,Q_{\lambda\varphi}(u)\geq-k^*\lambda^4\,\|u\|_2^2+T(u),
\]
for a form $T\in\cL$ and all $u\in H^2_0(\Omega)$.
\la{lemma6}
\end{lemma}
{\em Proof.} The fact that $\varphi\in\cE_{A,M}$ implies that $A(x,\nabla\varphi(x))\leq1$ for all $x\in\Omega$. Hence using Lemmas \ref{lem2}, \ref{lem:4} and \ref{lem3} we obtain
\bean
\re\,Q_{\lambda\varphi}(u)&=&-\int_\Omega k(x)A(x,\lambda\nabla\varphi)\,|u|^2\,dx+\Gamma_{\lf}(u)+T(u)\\
&\geq&-k^*\lambda^4\int_\Omega|u|^2\,dx+\Gamma_{\lf}(u)+T(u),
\eean
for some form $T\in\cL$ and all $u\in H^2_0(\Omega)$. Moreover
\bean
\Gamma_{\lf}(u)&=&\frac{1}{(2\pi)^{2}}\iiint_{\Omega\times\R^2\times\R^2}\Gamma(x, \, p_{x,\xi,\lambda\nabla\varphi},p_{x,\xi',\lambda\nabla\varphi})e^{i(\xi-\xi')\cdot x}\hat{u}(\xi)\overline{\hat{u}(\xi')}\,d\xi\,d\xi'\,dx\\
&=&\frac{1}{(2\pi)^{2}}\int_{\Omega}\Gamma\Big(x , \,\int_{\R^2}e^{i\xi\cdot x}\hat{u}(\xi)p_{x,\xi,\lambda\nabla\varphi}d\xi,\int_{\R^2}e^{i\xi'\cdot x}\hat{u}(\xi')p_{x,\xi',\lambda\nabla\varphi}d\xi'\Big)\,dx\\
&\geq&0,
\eean
by the positive definiteness of $\Gamma$; the result follows.$\hfill\Box$

We can now prove Theorem \ref{thm1}. We first consider the case where the coefficients of $H$ are Lipschitz continuous. For the general case
we shall then use the fact that Lemma \ref{lem:jde} is stable under $L^{\infty}$ perturbation of the coefficients.

{\bf\em Proof of Theorem \ref{thm1}} {\em Part 1.} We assume that the functions $\alpha(x),\beta(x),\gamma(x)$ are Lipschitz continuous. We claim that for any $\epsilon$ and $M$ positive there exists $c_{\epsilon,M}$ such that
\be
\re\,Q_{\lf}(u)\geq-\Big\{(k^*+\epsilon)\lambda^4+c_{\epsilon,M}(1+\lambda^3)\Big\}\|u\|_2^2.
\la{claim}
\ee
for all $\lambda>0$ and $\varphi\in\cE_{A,M}$.
To prove this we first note (cf. \cite[Lemma 7]{b2001}) that any form $T\in\cL$ satisfies
\[
|T(u)|\leq\epsilon Q(u)+c_{\epsilon}(1+\lambda^3)\,\|u\|_2^2,
\]
for all $\epsilon\in(0,1)$, $\lambda>0$ and $u\in H^2_0(\Omega)$. Hence Lemma \ref{lemma6} implies
\be
\la{gui1}
\re\,Q_{\lf}(u)\geq-\Big\{k^*\lambda^4+c_{\epsilon,M}(1+\lambda^3)\Big\}\|u\|_2^2-\epsilon Q(u).
\ee

Now, from (\ref{fil}) we have that for any $\psi\in\cE$ there holds
\[
\big| Q(u)-Q_{\psi}(u) \big|\leq\frac{1}{2}Q(u)+c\{\|\psi\|_{W^{2,\infty}}+\|\psi\|_{W^{2,\infty}}^4\}\|u\|_2^2.
\]
Taking $\psi=\lf$ where $\lambda>0$ and $\varphi\in\cE_{A,M}$ we thus obtain
\be
\big| Q(u)-Q_{\lf}(u) \big| \leq\frac{1}{2}Q(u)+c_M(\lambda+\lambda^4)\|u\|_2^2, 
\la{fil1}
\ee
Now, the coefficients of $\lambda^4$ in the expansion of $Q_{\lf}$ only involve first derivatives of $\varphi$. Since $|\nabla\varphi|\leq c$ for all $\varphi\in\cE_{A,M}$, (\ref{fil1}) can be improved to
\[
\big|Q(u)-Q_{\lf}(u) \big|\leq\frac{1}{2}Q(u)+\big\{c_M(\lambda+\lambda^3)+c\lambda^4\big\}\|u\|_2^2,
\]
which in turn implies
\be
Q(u)\leq\,  2{\rm Re}\,Q_{\lf}(u)+\big\{c_M(\lambda+\lambda^3)+c\lambda^4\big\}\|u\|_2^2.
\la{last}
\ee
Let $u\in H^2_0(\Omega)$ be given. If $\re Q_{\lf}(u)\geq0$ then (\ref{claim}) is obviously true. If not we then have from (\ref{gui1})
and (\ref{last})
\bean
\re\,Q_{\lf}(u)&\geq&-\Big\{k^*\lambda^4+c_{\epsilon}(1+\lambda^3)\Big\}\|u\|_2^2-2\epsilon\, \re Q_{\lf}(u)-\epsilon\big\{c_M(\lambda+\lambda^3)+c\lambda^4\big\}\|u\|_2^2\\
&\geq&-\Big\{(k^*+c\epsilon)\lambda^4+c_{\epsilon}(1+\lambda^3)+\epsilon\big\{c_M(\lambda+\lambda^3)+c\lambda^4\big\}\Big\}\|u\|_2^2,
\eean
and (\ref{claim}) again follows; hence the claim has been proved.

We complete the standard argument; Lemma \ref{lem:jde} and (\ref{claim}) imply
\[
|G(x,x',t)|<c_\epsilon t^{-1/2}\exp\Big\{\lambda[\varphi(x)-\varphi(x')]+(1+\epsilon)\big\{(k^*+\epsilon)\lambda^4+c_{\epsilon,M}(1+\lambda^3)\big\}\,t\Big\}, 
\]
for all $\epsilon\in(0,1)$. Optimizing over $\varphi\in\cE_{A,M}$ yields
\[
|G(x,x',t)|<c_\epsilon t^{-1/2}\exp\Big\{-\lambda d_M(x,x')+(1+\epsilon)\{(k^*+\epsilon)\lambda^4+c_{\epsilon,M}(1+\lambda^3)\}\,t\Big\} .
\]
Finally choosing
$
\lambda=[d_M(x,x')/(4k^*t)]^{1/3},
$
we have
\[
-\lambda d_M(x,x')+k^*\lambda^4t =-\sigma_*\frac{d_M(x,x')^{4/3}}{t^{1/3}},
\]
and (\ref{eq:mainthm}) follows.

{\em Part 2.} We now consider the general case where the functions $\alpha(x),\beta(x),\gamma(x)$ are not Lipschitz continuous.
Then there exist Lipschitz functions $\tilde\alpha(x),\tilde\beta(x),\tilde\gamma(x)$ such that (cf. (\ref{D}))
\[
\max\{\|\alpha-\tilde\alpha\|_{\infty},\;\|\beta-\tilde\beta\|_{\infty},\;\|\gamma-\tilde\gamma\|_{\infty}\}<2D.
\]
We assume that $D$ is small enough so that the corresponding operator $\tilde H$ is elliptic; we shall use a tilde to denote the various entities associated to $\tilde H$. 
Given $\varphi\in\cE_{\tilde{A},M}$ and $\lambda>0$ it follows from the first part of the proof that
\be
\la{gui5}
\re\,\tilde{Q}_{\lf}(u)\geq-\Big\{\tilde{k}^*\lambda^4+c_{\epsilon}(1+\lambda^3)\Big\}\|u\|_2^2-\epsilon Q(u).
\ee
Moreover it is easily seen that
\be
|k^*-\tilde{k^{*}}| \leq  cD\;,\qquad  \big| Q_{\lf}(u)-\tilde{Q}_{\lf}(u) \big| \leq cD\big\{Q(u)+\lambda^4\|u\|_2^2\big\}.
\la{gui6}
\ee
From (\ref{gui5}) and (\ref{gui6}) we obtain
\[
\re\,Q_{\lf}(u)\geq-\Big\{(k^*+cD)\lambda^4+c_{\epsilon}(1+\lambda^3)\Big\}\|u\|_2^2-\epsilon Q(u).
\]
As in Part 1, this leads to a Gaussian estimate involving the constant $\sigma_*-cD$ and the distance $\tilde{d}_M(x,x')$. To replace
$\tilde{d}_M(x,x')$ by $d_M(x,x')$ we note that there exists $c>0$ such that if $\varphi\in\cE_{\tilde{A},M}$ then $(1+cD)^{-1}\varphi\in \cE_{A,M}$. This implies that $\tilde{d}_M(x,x')\geq(1+cD)^{-1}d_{M}(x,x')$, which completes the  proof of the theorem.$\hfill\Box$

%%%%%%%%%%%%%%%%%%%%%%%%%%%%%%%%%%%%%%%%%%%%%%%%%%%%%%%%%%%%%%%%%%%%%%%%%%%%%%%%%%%%%%%%%%%%
%%%%%%%%%%%%%%%%%%%%%%%%%%%%%%%%%%%%%%%%%%%%%%%%%%%%%%%%%%%%%%%%%%%%%%%%%%%%%%%%%%%%%%%%%%%%
%%%%%%%%%%%%%%%%%%%%%%%%%%%%%%%%%%%%%%%%%%%%%%%%%%%%%%%%%%%%%%%%%%%%%%%%%%%%%%%%%%%%%%%%%%%%
%%%%%%%%%%%%%%%%%%%%%%%%%%%%%%%%%%%%%%%%%%%%%%%%%%%%%%%%%%%%%%%%%%%%%%%%%%%%%%%%%%%%%%%%%%%%

\section{Short time asymptotics}
\la{sec:asympt}

In this section we study the short time asymptotic behavior of the Green function $G(x,t)$ of the constant coefficient equation
\be
u_t=-\big( \partial_{x_1}^4+2\beta\partial_{x_1}^2\partial_{x_2}^2+\partial_{x_2}^4\big)u,\;\qquad x\in\R^2\;,\quad t>0.
\la{me}
\ee
(The slightly more general case where we have $\alpha\partial_{x_1}^4+2\beta\partial_{x_1}^2\partial_{x_2}^2+\gamma\partial_{x_2}^4$
is easily reduced to (\ref{me}).) The symbol of the elliptic operator is
\[
A(\xi) =\xi_1^4 +2\beta \xi_1^2\xi_2^2 +\xi_2^4 \; , \quad \xi \in\R^2,
\]
and it is strongly convex if and only if $0 < \beta < 3$.

Theorem \ref{thm2} below implies the sharpness of the constant $\sigma_*$ of Theorem \ref{thm1}, but it is interesting on its own.
As already mentioned, the behavior when $\beta\leq 0$ or $\beta\geq 3$ is qualitatively different from that of the case $0<\beta <3$
studied in \cite{ep}. The borderline cases $\beta=0,3$ are particularly interesting.

The Green's function for equation (\ref{me}) is given by
\be
G(x,t)=\frac{1}{(2\pi)^2}\int_{\R^2}e^{i\,\xi\cdot\,x-t\,A(\xi)}\,d\xi,\qquad x\in\R^2\, ,\; t>0.
\la{gf}
\ee
As already noted in the Introduction, we only consider specific points $x$: points lying on any coordinate axis when $\beta\geq 3$ and points lying on any main bisector when $\beta \leq 0$; due to symmetries this amounts to points of the form
$(s,0)$ and $(s,s)$ respectively. This choice is related to Lemma \ref{lem1}: in each of the two cases the respective point $\eta$ (i.e. $\eta =(s,0)$ or $\eta =(s,s)$) is a point for which there exists $\xi\in\R^2$ so that (\ref{eq1}) becomes an equality.
Moreover, for these points the explicit computation of the distance to the origin is possible; see also Subsection \ref{nss} below.

The main result of this section is the following
\begin{theorem}
For any $s>0$ the following asymptotic formulae are valid as $t\to 0+$:
\bean
& \ia & \mbox{ If $-1<\beta<0$ and $x=(s,s)$ then} \\
&& G(x ,t) \sim \frac{1}{3^{1/2}\cdot 4^{1/3}\, \pi} \frac{(1-\beta)^{1/6}}{(3-\beta)^{1/2}(1+\beta)^{1/6}}s^{-2/3} t^{-1/3}
\exp\Big(  -\frac{3}{ 4^{4/3}}\big(  \frac{1+\beta}{1-\beta}\big)^{1/3}\frac{s^{4/3}}{t^{1/3}} \Big) \\[0.2cm]
& \ib & \mbox{ If $\beta=0$ and $x=(s,s)$ then} \\
&& G(x ,t) \sim \frac{1}{3 \cdot 4^{1/3}\, \pi} s^{-2/3} t^{-1/3}
\exp\Big(  -\frac{3}{ 4^{4/3}} \frac{s^{4/3}}{t^{1/3}} \Big) \cdot 
\Big( 1+ \cos\Big[ \frac{3\sqrt{3}}{ 4^{4/3}} \frac{s^{4/3}}{t^{1/3}} -\frac{\pi}{3}  \Big] \Big)\\[0.2cm]
& \ic & \mbox{ If $\beta=3$ and $x=(s,0)$ then} \\
&&  G (x,t )\sim 
\frac{1}{3 \cdot 4^{1/3} \, \pi} s^{-2/3}t^{-1/3} \exp\Big( -\frac{3}{8\cdot 4^{1/3}}\frac{s^{4/3}}{t^{1/3}}\Big) \cdot
\Big( 1+ \cos\Big[ \frac{3\sqrt{3}}{8\cdot 4^{1/3}}\frac{s^{4/3}}{t^{1/3}} -\frac{\pi}{3}  \Big] \Big) \\[0.2cm]
& \id & \mbox{ If $\beta>3$ and $x=(s,0)$ then} \\
&&  G(x,t) \sim  \frac{1}{ 2^{7/6} \cdot 3^{1/2}\, \pi} \beta^{-1/2}(\beta^2-1)^{1/6}s^{-2/3} t^{-1/3} 
\exp\Big( -\frac{3}{4^{4/3}} (\beta^2-1)^{-1/3}\frac{s^{4/3}}{t^{1/3}} \Big)  
\eean
\la{thm2}
\end{theorem}

{\bf Remark.} Clearly that the notation $F(\lambda)\sim G(\lambda)$ cannot have here the usual meaning $F(\lambda) =G(\lambda)( 1+o(\lambda))$, as the function $G$ takes also the value zero.
By looking at the proof below it becomes clear that the actual meaning of
\[
F(\lambda ) \sim e^{A \lambda} \lambda^D \big[1+ \cos(B\lambda +C) \big]   \; , \quad \mbox{ as }\lambda \to +\infty,
\]
is that
\[
F(\lambda ) =  e^{A \lambda}  \lambda^D  \Big( [1+ \cos(B\lambda +C) ]+o(1) \Big) \; , \quad \mbox{ as }\lambda \to +\infty.
\]

\subsection{Some comments on the distance $d(x,x')$}
\la{nss}

Before proceeding with the proof of Theorem \ref{thm2} we make some comments on the distance $d(x,x')$ defined by (\ref{fd}). First, we recall that a Finsler metric on a domain $\Omega\subset\R^n$ is map $p:\Omega \times \R^n \to \R_+$ whose regularity with respect to $x\in\Omega$ may vary and which has the following properties
\bean
&\ia  & p(x,\xi)=0  \mbox{ if and only if }\xi =0  \\
&\ib &  p(x,\lambda\xi)= |\lambda|  p(x,\xi) \; , \quad \lambda\in\R  \\
& \ic & p(x,\xi) \mbox{ is convex with respect to }\xi
\eean
Given a Finsler metric on $\Omega$ the dual metric $p_*$ is defined by
\be
p_*(x,\eta) =  \max_{\xi \neq 0} \frac{ \eta\cdot \xi}{ p(x,\xi)} \; , \qquad x\in\Omega \; , \;\;\; \eta\in\R^2 \, .
\la{pstar}
\ee
This is also a Finsler metric and there holds $p_{**}=p$. Having a Finsler metric one can define lengths of paths and hence the (geodesic) distance between points.

We now return to our specific case. The map
\[
 p(x, \xi) =A(x, \xi)^{1/4}  \; , \quad  x\in\Omega \, , \;\; \xi\in\R^2,
\]
satisfies properties (i) and (ii) above but not (iii). Nevertheless the dual metric $p_*$ can still be defined by (\ref{pstar}). Since it is convex (being the supremum of linear functions) it is a Finsler metric. Clearly $p_{**}(x,\xi)$ does not coincide with $p(x,\xi)$ in this case. Actually, there holds $p_{**}(x,\xi)\leq p(x,\xi)$; indeed it may be seen that the set $\{ \xi : p_{**}(x,\xi) \leq 1 \}$ is precisely the convex hull
of the set $\{ \xi : p(x,\xi) \leq 1\}$ .

Now, the (geodesic) Finsler distance $d_{p_*}(x,x')$ induced by $p_*$ satisfies \cite[Lemma 1.3]{agmon}
\[
d_{p_*}(x,x') =   \sup \{ \varphi(x') -\varphi(x) \; : \;  \varphi\in{\rm Lip}(\Omega) , \;\; p_{**}( y ,\nabla\varphi(y))\leq 1 \, , \;\;   \mbox{ a.e. } y\in\Omega\}.
\]
Since $p_{**}\leq p$ this implies $d(x,x')\leq d_{p_*}(x,x')$. We shall now see that this does not spoil the sharpness of the constant $\sigma_*$ of Theorem \ref{thm1}.

Let us restrict from now on our attention to the constant coefficient case.
By translation invariance we have $d(x,x') =d_0(x-x')$ where
\be
d_0(x) = \sup \{ \varphi(x) \; : \;  \varphi\in{\rm Lip}(\R^2) , \;\;   \varphi(0)=0  \;  , \;\;   A( y ,\nabla\varphi(y))\leq 1 \, , \;\;   \mbox{ a.e. } y\in\R^2\}.
\la{fd0}
\ee
We then have
\be
d_0(x) =p_*(x) \; , \qquad x\in\R^2.
\la{lolo}
\ee
Indeed, given a function $ \varphi$ as in (\ref{fd0}) we have
\[
\varphi(x) =\int_0^1 \frac{d}{dt} \varphi(tx) dt = \int_0^1 \nabla\varphi(tx)\cdot x \,  dt \leq  \int_0^1 A\big(\nabla\varphi(tx)\big)^{1/4} p_*(x) \,  dt  \leq p_*(x) \,  ,
\]
hence $d_0(x)\leq p_*(x)$. For the converse, let $\xi\in\R^2\setminus\{0\}$ be given. The function
\[
\varphi(y) = \frac{\xi\cdot y}{ A(\xi)^{1/4} } \; \; , \qquad y\in\R^2,
\]
then satisfies $A(\nabla \varphi(y))=1$ and therefore can be used as a test function in (\ref{fd0}).
Hence
\[
\frac{\xi\cdot x}{ A(\xi)^{1/4}}  =\varphi(x) \leq d_0(x),
\]
and maximizing over $\xi$ yields $p_*(x)\leq d_0(x)$.

Now, it is immediate from (\ref{pstar}) that
\be
 p_*(x)  \geq  \frac{|x|^2}{p(x)}   \; , \qquad  x\in\R^2\setminus\{0\}.
\la{di}
\ee
We shall need to identify the points $x\in\R^2$ for which (\ref{di}) becomes an  equality. By homogeneity it is enough to consider points of unit Euclidean length.
Let us write $e_{\phi}=(\cos\phi, \sin\phi)$. We are then seeking directions $\phi$ for which
\[
 \frac{1}{ A(e_{\phi})^{1/4} }  =\max_{\theta\in\R} \frac{ e_{\phi}\cdot e_{\theta}}{ A(e_{\theta})^{1/4} } .
\]
So let $\phi\in [0,2\pi]$ be fixed and set
\[
g(\theta) = \frac{ e_{\phi}\cdot e_{\theta}}{ A(e_{\theta})^{1/4} }=  \frac{\cos(\phi-\theta)}{ \Big(   1 +\frac{\beta -1}{2}\sin^2 2\theta \Big)^{1/4}} .
\]
Then
\[
g'(\theta) =   A(e_{\theta})^{-1/4}  \sin (\theta -\phi)  - \frac{\beta-1}{4} A(e_{\theta})^{-5/4} \cos(\theta-\phi) \,  \sin 4\theta \, .
\]
It follows that $g'(\phi)=0$ if and only if $ \sin 4\phi =0$, i.e. if and only if $\phi$ is an integer multiple of $\pi/4$. This corresponds exactly to the points considered in Theorem \ref{thm2} and hence for these points inequality (\ref{di})
holds as an equality. In particular, recalling (\ref{lolo}) we have
\[
d_0(s,0) = s
\]
and
\[
d_0(s,s) = \frac{2s^2}{ (2s^4 +2\beta s^4)^{1/4}} =  2^{3/4} (1+\beta)^{-1/4}s \, .
\]

Suppose now that  $\beta\geq 3$. We then have  (cf. (\ref{sx})) $\sigma = 3 \cdot 4^{-4/3} (\beta^2-1)^{-1/3}$ (recall that $Q=\beta$ for equation (\ref{me})).
This proves the sharpness of the constant $\sigma_*$ in Theorem \ref{thm1} in the regime $Q\geq 3$.

Similarly, if $-1< \beta \leq 0$ then (cf. (\ref{sx}))
\[
\sigma = \frac{3}{2 \cdot 4^{4/3}} \frac{(1+\beta)^{2/3}}{(1-\beta)^{1/3}}
\]
and therefore
\[
\exp\Big(  -\frac{3}{4^{1/3}}\big(  \frac{1+\beta}{1-\beta}\big)^{1/3}\frac{s^{4/3}}{t^{1/3}} \Big) =\exp\Big(-\sigma_* \frac{d_0(s,s)^{4/3}}{t^{1/3}}\Big).
\]
Hence the constant $\sigma_*$ is sharp also in the regime $-1< Q\leq 0$.

\subsection{Proof of Theorem \ref{thm2}}

Changing variables in (\ref{gf}) by $\xi=(4t)^{-1/3} \eta$ we obtain
\be
G(x,t)=\frac{1}{(2\pi)^2 (4t)^{2/3}}\, F\Big( \frac{1}{(4t)^{1/3}}\Big),
\la{cv}
\ee
where
\[
F(\lambda)=\int_{\R^2}e^{ \lambda\big( i\, x\cdot\,\xi-\frac{1}{4}A(\xi) \big)}\,d\xi .
\]
To find the asymptotic behavior of $F(\lambda)$ as $\lambda\to +\infty$ we shall use the method of steepest descent. So we shall consider the complex analytic function of two variables, $z=(z_1,z_2)$,
\[
\phi(z)=i\,x\cdot z-\frac{1}{4}A(z)  =i(x_1z_1 + x_2z_2) -\frac{1}{4}(z_1^4+2\beta z_1^2z_2^2 +z_2^4),
\]
and shall use Cauchy's theorem for functions of two variables to suitably deform $\R^2\subset\C^2$ to some other surface in $\C^2$ that will contain the saddle points of $\phi$ that actually contribute to the asymptotic behavior of $F(\lambda)$.
For our purposes it is enough to consider deformations that are parallel transports by a point in $i\R^2$. Indeed, it easily follows by Cauchy's theorem that for any $\eta_0\in\R^2$ we have
\[
F(\lambda) =\int_{\R^2 +i\eta_0}  e^{\lambda \phi(z)}dz.
\]
The main issue is to identify the relevant saddle points and (hence) the vector $\eta_0$. What is of importance here is the real part of
$\re \phi(z)$ -- also called the height of $\phi(z)$.
The relevant saddle points are not necessarily those of the largest height, but rather, they are those for which there exists a deformation such that the largest height on it is attained at those points.

Concerning the notation, we  shall write each $z\in\C^2$ as $z=(z_1,z_2)$ but also as
$z =\xi + i\eta$ with $\xi =(\xi_1,\xi_2)\in\R^2$ and $\eta =(\eta_1,\eta_2)\in\R^2$. Finally we note that it is enough to prove the asymptotic formulae
in case $s=1$, since the general case then follows from the relation
\[
G(s\, x, t) =\frac{1}{s^2} \, G(x, \frac{t}{s^4}) \; , \qquad t,s>0 \, , \;\; x\in\R^2.
\]

%%%%%%%%%%%%%%%%%%%%%%%%%%%%%%%%%%%%%%%%%%%%%%%%%%%

\subsubsection{The case $-1 < \beta\leq 0$}

In this case we have $x=(1,1)$. Two saddle points that are relevant are the points
\[
z_0^{\pm} =  \pm \xi_0  +i\eta_0 
\]
where
\[
\xi_0 =   \frac{1}{2}    \frac{(3-\beta)^{1/2}}{(1+\beta)^{1/6}(1-\beta)^{1/3}} (1,-1)  \; , 
\qquad\quad
\eta_0 = \frac{1}{2} \Big(  \frac{1+\beta}{1-\beta} \Big)^{1/3}  (1,1) \, .
\]
We deform $\R^2$ by $i\eta_0$ and have
\[
F(\lambda) =\int_{\R^2+i\eta_0} e^{\lambda\phi(z)}dz \, .
\]

{\bf Case 1. $-1<\beta<0$.} In this case the saddle points that contribute are precisely the points $z_0^{\pm}$. 
We claim that
\be
\re \phi(z)  \leq \re \phi(z_0^+)  \; , \qquad z\in\R^2+i\eta_0,
\la{fi1}
\ee
with equality exactly at the points $z_0^{\pm}$. To prove (\ref{fi1}) we note that it is equivalently written as
\[
\re A(z)  \geq \re A(z_0^+)  =    -\Big(  \frac {1+\beta}{1-\beta} \Big)^{1/3}    ,
\]
so it is enough to establish that
\[
\re A(\xi +i\eta_0) +  \Big(  \frac {1+\beta}{1-\beta} \Big)^{1/3}  \geq 0 \; , \qquad \xi\in\R^2.
\]
This is indeed true, as a direct computation shows that
\bea
\re A(\xi +i\eta_0) + \Big(  \frac {1+\beta}{1-\beta} \Big)^{1/3} & =& -\beta(\xi_1^2 -\xi_2^2)^2  +
(\beta +1) \bigg[   \Big(  \xi_1^2 - \frac{3-\beta}{4(1+\beta)^{1/3}(1-\beta)^{2/3}}\Big)^2  \nonumber \\
&& +   \Big(  \xi_2^2 - \frac{3-\beta}{4(1+\beta)^{1/3}(1-\beta)^{2/3}}\Big)^2  \bigg]   -\beta \Big( \frac{ 1+\beta}{1-\beta} \Big)^{2/3} (\xi_1+\xi_2)^2  \nonumber \\
&\geq & 0,
 \la{la2}
\eea
[This is a scaled version of (\ref{relation1}) for $\eta =(1,1)$.]
Clearly equality holds only for the points $\xi=\pm \xi_0$, and these correspond to the points $z_0^{\pm}$.
Hence the claim has been proved.

This implies (see \cite[Criterion 1, page 15]{ep}) that the points $z_0^{\pm}$ are precisely those that contribute to the asymptotic behavior of $F(\lambda)$ as $\lambda\to +\infty$. Now, it is easy to see that the for any $\delta>0$ the integrals
\[
 \int_{D(z_0^{\pm} ,\delta) +i\eta_0}e^{\lambda \phi(z)}dz
\]
are complex conjugate of each other, hence the total contribution of the these two points is equal to twice the real part of the contribution of $z_0^{+}$. Since these saddle
points are non-degenerate, the contribution of $z_0^+$ is given by the formula (see \cite[equation (3.6)]{ep} or \cite[equation (1.61)]{f})
\[
{\rm contr}(z_0^+) =  \frac{2\pi}{\lambda} \big( {\rm det} (\phi_{z_iz_j})\big|_{z=z_0^+} \big) ^{-1/2} e^{\lambda \phi(z_0^+)}.
\]
We have
\[
\phi(z_0^+) =-\frac{3}{4}\Big( \frac{1+\beta}{1-\beta}\Big)^{1/3}   \;  , \qquad 
 {\rm det} (\phi_{z_iz_j})|_{z=z_0^+} = \frac{3(3-\beta)(1+\beta)^{1/3}}{(1-\beta)^{1/3}},
\]
hence combining the above we conclude that
\be
F(\lambda)\sim \frac{4\pi}{\lambda} \frac{(1-\beta)^{1/6}}{\sqrt{3} (3-\beta)^{1/2}(1+\beta)^{1/6}}
\exp\Big(  - \frac{3}{4}\big( \frac{1+\beta}{1-\beta} \big)^{1/3} \lambda \Big) , \qquad
\mbox{ as }\lambda \to +\infty.
\la{888}
\ee
Recalling (\ref{cv}) concludes the proof in this case.

{\bf Case 2. $\beta=0$.} In this case $G(x,t)$ is the square of an one-dimensional integral; we prefer however to use the two-dimensional
approach because the setting is already prepared, but also because we believe that this conveys better the essential issues involved.

Relation (\ref{la2}) is also valid for $\beta=0$ in which case it is written
\[
\re A(\xi +i\eta_0) + 1  =
 \big(  \xi_1^2 - \frac{3}{4}\big)^2  +   \big(  \xi_2^2 - \frac{3}{4}\big)^2   \geq  0 \, .
\]
The points $z_0^{\pm}$ considered above are saddle points also for $\beta=0$. The same computations as above are valid hence their contribution is (cf. (\ref{888}))
\[
{\rm contr}(z_0^+)+{\rm contr}(z_0^-) = \frac{4\pi}{3\lambda} \exp\big( -\frac{3}{4}\lambda\big).
\]
However in this case there are two more saddle points of $\phi$ that lie on $\R^2+i\eta_0$ and that must be considered, namely the points
\[
z_*^{\pm}  =\pm\frac{\sqrt{3}}{2}(1,1)   + i\eta_0 \, ,
\]
For these points we find
\[
\phi(z_*^{\pm}) = -\frac{3}{4}\pm\frac{3\sqrt{3}}{4} i   \;  , \qquad  {\rm det} (\phi_{z_iz_j})|_{z=z_*^{\pm}} = 9e^{2\pi i/3},
\]
and thus obtain the contribution
\[
{\rm contr}(z_*^+)+{\rm contr}(z_*^-) = \frac{4\pi}{3\lambda} \exp\big( -\frac{3}{4}\lambda\big)
\cos\Big( \frac{3\sqrt{3}}{4} \lambda   -\frac{\pi}{3}\Big).
\]
Adding the contributions we arrive at
\be
 F(\lambda) \sim \frac{4\pi}{3\lambda} \exp\big( -\frac{3}{4}\lambda\big)  \Big( 1+
\cos\Big( \frac{3\sqrt{3}}{4} \lambda   -\frac{\pi}{3}\Big) \Big) ,
\la{last1}
\ee
which concludes the proof by means of (\ref{cv}). $\hfill\Box$

\subsubsection{The case $\beta\geq 3$}

In this case we have $x=(1,0)$. Two saddle points that are relevant in this case are the points
\[
z_0^{\pm}=(\beta^2-1)^{-1/3}[(0,\pm\sqrt{\beta})+i(1,0)] = \pm \xi_0 +i\eta_0 \, .
\]
As before, we have
\[
F(\lambda) =\int_{\R^2+i\eta_0} e^{\lambda\phi(z)}dz \, .
\]

{\bf Case 1. $\beta>3$.} In this case the relevant saddle points are precisely the points $z_0^{\pm}$. This will follow if we prove that
\be
\re \phi(z)  \leq \re \phi(z_0^+)  \; , \qquad z\in\R^2+i\eta_0,
\la{fi}
\ee
with equality exactly at the points $z_0^{\pm}$. To prove (\ref{fi}) we note that it is equivalently written as
\[
\re A(z)  \geq \re A(z_0^+)  = -(\beta^2-1)^{-1/3},
\]
so it is enough to establish that
\[
\re A(\xi +i\eta_0) + (\beta^2-1)^{-1/3}  \geq 0 \; , \qquad \xi\in\R^2.
\]
This is indeed true, as a direct computation shows that
\bea
&& \hspace{-1.5cm}\re A(\xi +i\eta_0) + (\beta^2-1)^{-1/3} \nonumber \\
&& =\Big(\xi_1^2  + \xi_2^2-\beta(\beta^2-1)^{-2/3} \Big)^2   +2(\beta-1)\xi_1^2\xi_2^2  + 2(\beta-3)(\beta^2-1)^{-2/3} \xi_1^2 \geq 0 .
\la{la}
\eea
[This is a scaled version of (\ref{relation3}) for $\eta =(1,0)$.] Equality holds only for the points
\[
\xi_0^{\pm} =\pm  \big(0, \sqrt{\beta}(\beta^2-1)^{-1/3} \big) 
\]
which correspond to the points $z_0^{\pm}$.

The two contributions are again complex conjugate of each other. We use again the relation
\[
{\rm contr}(z_0^+) =  \frac{2\pi}{\lambda} \big( {\rm det} (\phi_{z_iz_j})\big|_{z=z_0^+} \big) ^{-1/2} e^{\lambda \phi(z_0^+)} ,
\]
and since 
\[
\phi(z_0^+) = -\frac{3}{4}(\beta^2 -1)^{-1/3} \; , \qquad\quad  {\rm det} (\phi_{z_iz_j})\big|_{z=z_0^+}  =6\beta(\beta^2-1)^{-1/3},
\]
combining the above we obtain
\be
F(\lambda)  \sim   \frac{4\pi}{\lambda}  (6\beta)^{-1/2} (\beta^2-1)^{1/6}  \exp\Big(-\frac{3}{4}(\beta^2-1)^{-1/3} \lambda\Big) \; , \qquad \mbox{ as }\lambda\to +\infty .
\la{that}
\ee
The proof is concluded by using (\ref{cv}).

{\bf Case 2. $ \beta=3$.} Inequality (\ref{la}) is also true for $\beta=3$, in which case it takes the form.
\[
\re A(\xi +i\eta_0) + \frac{1}{2} 
 =\big(\xi_1^2  + \xi_2^2- \frac{3}{4} \big)^2   +4 \xi_1^2\xi_2^2   \geq 0
\]
In this case equality holds not only at the points $\xi_0^{\pm}$ but also at the points
\[
\xi_*^{\pm} =\pm  \big( \frac{\sqrt{3}}{2} \, ,  \, 0 \big) .
\]
The corresponding points in $\C^2$ are the points
\[
 z_*^{\pm} =\xi_*^{\pm} +i\eta_0.
\]
As before, the combined contribution of the points
$ z_*^{\pm}$ is twice the real part of the contribution of $ z_*^{+}$. We find
\[
\phi(z_*^+) =-\frac{3}{8} +\frac{3\sqrt{3}}{8} i   \;  , \qquad  {\rm det} (\phi_{z_iz_j})|_{z=z_*^+} =9 e^{ 2 \pi i/3},
\]
hence using the same formula as above we obtain
\[
{\rm contr}(z_*^+)+{\rm contr}(z_*^-) = \frac{4\pi}{3\lambda} \exp \big( -\frac{3}{8} \lambda \big) 
\cos\big( \frac{3\sqrt{3}}{8}\lambda +\frac{\pi}{3}  \big).
\]
The contribution of the first two points $z_0^{\pm}$ is given by (\ref{that}) (for $\beta=3$); adding the two contributions we conclude that
\be
F(\lambda) \sim \frac{4\pi}{3\lambda} \exp \big( -\frac{3}{8} \lambda \big) 
\Big( 1+ \cos\big[ \frac{3\sqrt{3}}{8}\lambda -\frac{\pi}{3}  \big] \Big) \;\; , \quad \mbox{ as }\lambda \to +\infty.
\la{that1}
\ee
The proof is concluded by recalling (\ref{cv}).  $\hfill\Box$

The estimates (\ref{888}), (\ref{last1}), (\ref{that}) and (\ref{that1}) obtained in the proof above all have the form $F(\lambda)\sim G(\lambda)$ for some explicitly given function $G(\lambda)$.
In each of the diagrams below we have plotted the numerically computed graph of $F(\lambda)e^{\sigma \lambda}$
(blue, dashed) against the function $G(\lambda)e^{\sigma \lambda}$ (red, continuous),
where $\sigma$ is the positive constant in the exponential term of
$G(\lambda)$. We note that in the case $\beta=4$ the convergence is slower, but more detailed computations are in line with the difference being
of order $O(1/\lambda^2)$.

\begin{figure*}[h]
\centering
\begin{minipage}{.5\textwidth}
  \centering
  \includegraphics[width=1\linewidth]{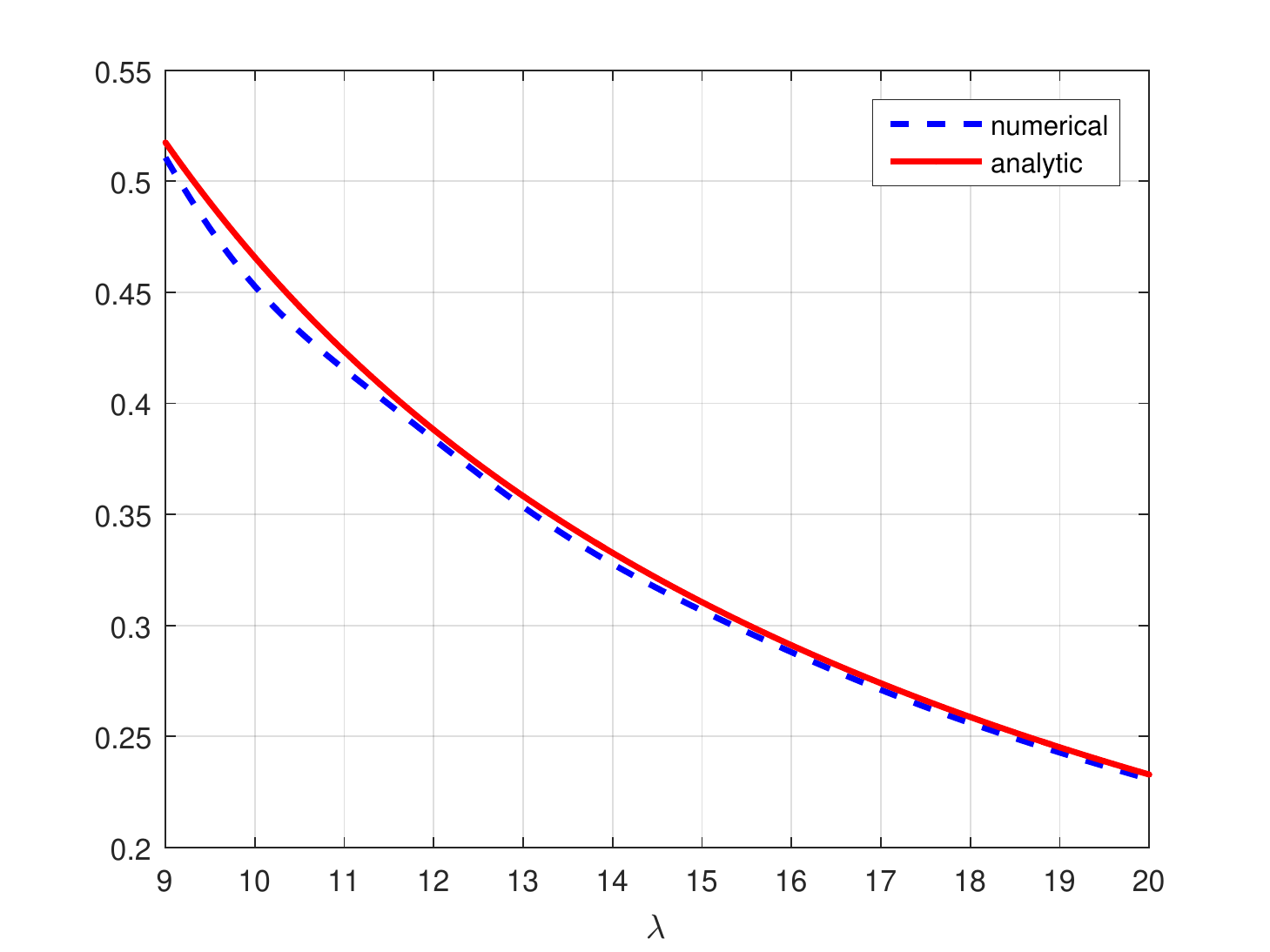}
  \caption{$\beta=-0.5$, $x=(1,1)$}
  %\label{fig:test1}
\end{minipage}%
\begin{minipage}{.5\textwidth}
  \centering
  \includegraphics[width=1\linewidth]{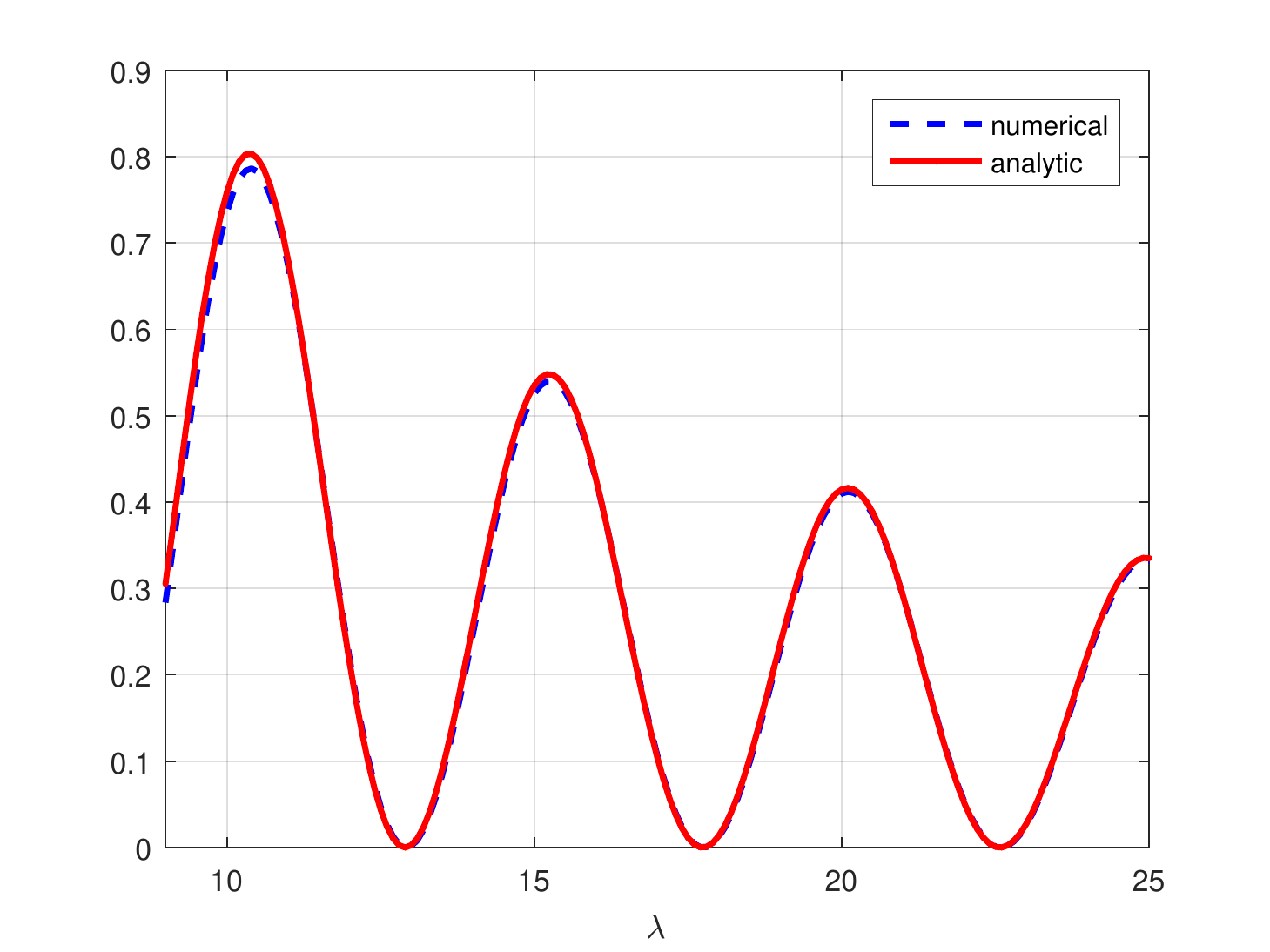}
  \caption{$\beta=0$, $x=(1,1)$}
  %\label{fig:test2}
\end{minipage}
\end{figure*}

\begin{figure*}[h]
\centering
\begin{minipage}{.5\textwidth}
  \centering
  \includegraphics[width=1\linewidth]{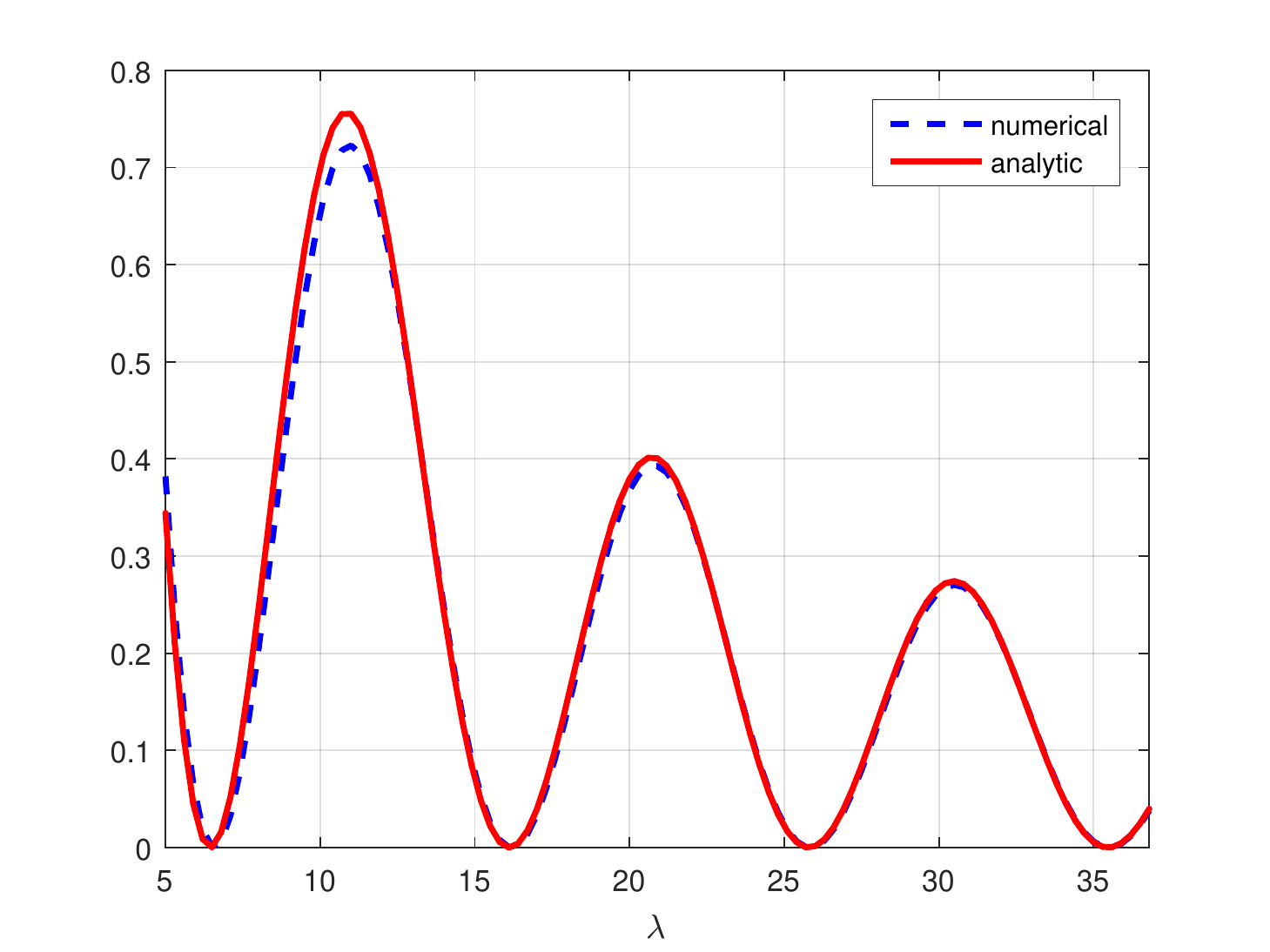}
  \caption{$\beta=3$, $x=(1,0)$}
  %\label{fig:test1}
\end{minipage}%
\begin{minipage}{.5\textwidth}
  \centering
  \includegraphics[width=1\linewidth]{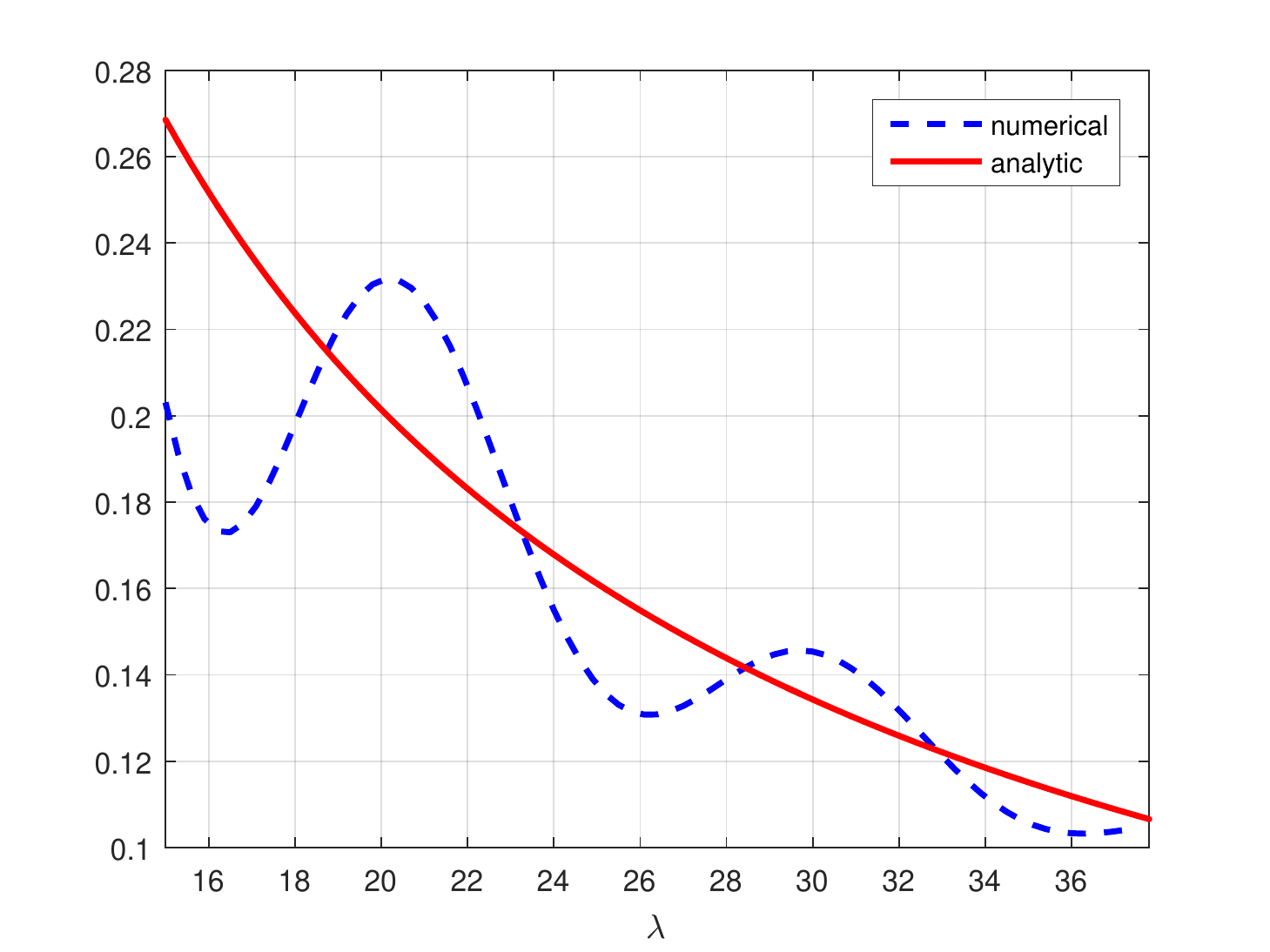}
  \caption{$\beta=4$, $x=(1,0)$}
  %\label{fig:test2}
\end{minipage}
\end{figure*}

%%%%%%%%%%%%%%%%%%%%%%%%%%%%%%%%%%%%%%%%%%%%%%%%%%%%%%%%%

%%%%%%%%%%%%%%%%%%%%%%%%%%%%%%%%%%%%%%%%%%%%%%%%%%%%%%%%%

\newpage

{\bf \Large Appendix}

In this appendix we present the proof of Evgrafov and Postnikov \cite{ep} for the asymptotic behavior of $G(x,t)$ in the strongly convex case $0<\beta< 3$
and for any $x\in\R^2$, $x\neq 0$.
We note that the article \cite{ep} deals with the general equation of an an operator of order $2m$ acting in $\R^d$.

To find the contributing saddle points we first note that by the strict convexity of the symbol there exists a unique $q=q(x)\in\R^2$ such that
\[
\frac{1}{4}\nabla A(q) =x.
\]
Then a point $z=\alpha q\in \C^2$ is a critical point for $\phi$ if and only if
$\alpha^3 =i$. We shall use two of these points, the points
\[
z_*^{\pm} = \big( \pm \frac{\sqrt{3}}{2} +\frac{1}{2} i   \big) q   =: \pm \xi_0 +i\eta_0 \, .
\]
As in the proof of Theorem \ref{thm2}  we change domain of integration from $\R^2$ to $\R^2+i\eta_0$ and the crucial property is that
\[
\re \phi(z) \leq \re \phi(z_0) \; , \qquad \mbox{ for all }z\in\R^2 +i\eta_0.
\]
with equality only at the points $z_*^{\pm}$.
To prove this we note that it is equivalent to
\[
\re A(\xi +i\eta_0 ) \geq \re A(\xi_0 + i\eta_0)   \, , \qquad \xi\in\R^2,
\]
with equality only for $\xi=\pm \xi_0$. 
With $\alpha$ as above, i.e. $\alpha =\pm \sqrt{3}/2 +i/2$, we compute
\[
\re A(\xi_0 + i\eta_0) = \re A(\alpha q)  = \re (\alpha^4) A(q) = -\frac{1}{2}A(q) =-8A(\eta_0),
\]
so we need to prove that
\[
\re A(\xi +i\eta_0 ) + 8A(\eta_0)\geq 0  , \qquad \xi\in\R^2,
\]
with equality at $\xi=\pm\xi_0$.
This is indeed true since for any $\eta=(\eta_1,\eta_2)$ we have
\bean
{\rm Re} \, A(\xi+i \eta)+8 A(\eta)&=&\frac{\beta}{3}\Big\{ (\xi_1^2-3\eta_1^2) + (\xi_2^2-3\eta_2^2)  \Big\}^2
+\frac{4\beta}{3}  (\xi_1\xi_2-3\eta_1\eta_2)^2  \nonumber \\
 &&  +\frac{3-\beta}{3}\Big\{    (\xi_1^2-3\eta_1^2)^2+ (\xi_2^2-3\eta_2^2)^2\Big\} \\
&\geq & 0.
\eean
Hence the asymptotic behavior will indeed result precisely from the points $z_*^{\pm}$. To compute it we first note that
\[
 \phi(z_*^{\pm}) = \frac{3}{4} e^{\pm \frac{2\pi i}{3}}A(q) =-\frac{3}{8}A(q) \pm \frac{3\sqrt{3}}{8}A(q)\, i.
\]
We also have
\[
{\rm det} (\phi_{z_iz_j})|_{z=z_*^+} = h(x)^{\frac{4}{3}}  e^{\frac{2\pi i}{3}}
\]
where the function $h$ is positively homogeneous of degree one. Hence
\[
{\rm contr}(z_*^+) =\frac{2\pi}{\lambda} h(x)^{-\frac{2}{3}} e^{-\frac{\pi i}{3}} 
\exp\Big(-\frac{3}{8}A(q)\lambda\Big) \exp\Big(\frac{3\sqrt{3}}{8}A(q)\, i\Big).
\]
The contribution of $z_*^-$ is the complex conjugate of that of $z_*^+$ and adding the two contributions we obtain that
\be
F(\lambda) \sim  \frac{4\pi}{\lambda} h(x)^{-2/3}
\exp\Big(-\frac{3}{8}A(q)\lambda\Big) \cos\Big(\frac{3\sqrt{3}}{8}A(q) -\frac{\pi}{3} \Big)
\la{finfin}
\ee
We claim that $A(q)=d_0(x)^{4/3}$. Indeed by (\ref{lolo}) we have
\[
d_0(x) =p_*(x) 
= \sup_{\xi}\frac{x\cdot \xi}{A(\xi)^{1/4}} 
\geq \frac{x\cdot q}{A(q)^{1/4}} 
= \frac{\frac{1}{4} \nabla A(q)\cdot q}{A(q)^{1/4}}
=A(q)^{3/4}
\]
The reverse inequality follows by noting that the supremum is attained at $\xi=q$. 

Substituting $A(q)=d_0(x)^{4/3}$ in (\ref{finfin}) and using (\ref{cv}) we finally conclude that as $t\to 0+$.
\[
G(x,t) \sim \frac{2^{1/3}}{\pi} 
h(x)^{-2/3} t^{-1/3}\exp\Big(-\frac{3}{8\cdot 4^{1/3}} \frac{d(x)^{\frac{4}{3}}}{t^{\frac{1}{3}}} \Big) 
\cos\Big(\frac{3\sqrt{3}}{8\cdot 4^{1/3}} \frac{d_0(x)^{\frac{4}{3}}}{t^{\frac{1}{3}}}  -\frac{\pi}{3}\Big) .
\]

\

{\bf Acknowledgment.} We thank Leonid Parnovski for useful suggestions and Gregory Kounadis for helping us with Matlab. We also thank the referee for crucial comments which led to a substantial improvement of
Section \ref{sec:asympt}.

\end{document}